\documentclass[english]{extarticle}
\usepackage[utf8]{inputenc}
\usepackage{geometry}
\usepackage{color}
\usepackage{float}
\usepackage{amsmath}
\usepackage{amsthm}
\usepackage{amssymb}
\usepackage{stmaryrd}
\usepackage{esint}
\PassOptionsToPackage{normalem}{ulem}
\usepackage{ulem}
\usepackage{cite}
\usepackage{graphicx}
\makeatletter

\floatstyle{ruled}
\newfloat{algorithm}{tbp}{loa}
\providecommand{\algorithmname}{Algorithm}
\floatname{algorithm}{\protect\algorithmname}

\theoremstyle{plain}
\newtheorem{thm}{\protect\theoremname}
\theoremstyle{plain}

\theoremstyle{definition}

\theoremstyle{remark}
\newtheorem{rem}[thm]{\protect\remarkname}
\theoremstyle{plain}

\theoremstyle{plain}

\ifx\proof\undefined
\newenvironment{proof}[1][\protect\proofname]{\par
\normalfont\topsep6\p@\@plus6\p@\relax
\trivlist
\itemindent\parindent
\item[\hskip\labelsep\scshape #1]\ignorespaces
}{%
\endtrivlist\@endpefalse
}
\providecommand{\proofname}{Proof}
\fi
\theoremstyle{remark}

\makeatother

\usepackage{babel}
\providecommand{\claimname}{Claim}
\providecommand{\definitionname}{Definition}
\providecommand{\lemmaname}{Lemma}
\providecommand{\propositionname}{Proposition}
\providecommand{\remarkname}{Remark}
\providecommand{\theoremname}{Theorem}
\providecommand{\corollaryname}{Corollary}
\newcommand{\x}{\textbf{x}}

\newcommand{\z}{\textbf{z}}
\newcommand{\w}{\textbf{w}}
\newcommand{\g}{\textbf{g}}
\newcommand{\f}{\textbf{f}}
\newcommand{\LL}{\textbf{L}}
\newcommand{\PH}{\boldsymbol{\Phi}}

\newcommand{\V}{\mathcal{V}}
\newcommand{\N}{\mathcal{N}}

\usepackage{latexsym}
\begin{document}

\title{Frank-Wolfe variants for minimization of a sum of functions }

\author{Suhail M. Shah}
\date{September 7, 2020}
\maketitle

\begin{abstract}
We propose several variants of the Frank-Wolfe algorithm to minimize a sum of functions. The main proposed algorithm is inspired from the dual averaging scheme of Nesterov adapted for Frank Wolfe in a stochastic setting. A distributed version of this scheme is also suggested. Additionally, we also propose a Frank-Wolfe variant  based on incremental gradient techniques. The convergence rates for all the proposed algorithms are established. The performance is studied on least squares regression and multinomial classification.  
\end{abstract}

\section{Introduction}
In this paper, we consider the optimization problem
\begin{equation}\label{prob1}
\min_{x \in \Omega} \big\{ f(x) := \mathbb{E}_{\xi}\, [f(x, \xi) ]  = \int f(x,\omega)\, dP(\omega) \big\},
\end{equation}
where $P$ is the probability law of the random variable $\xi$, $x \in \mathbb{R}^d$, $f\,:\, \mathbb{R}^d \to \mathbb{R}$ and $\Omega \subset \mathbb{R}^d$ is compact convex set. We also consider a closely related  problem to the above, where the objective function can be expressed as a sum of functions :
\begin{equation}\label{prob}
\min_{x \in \Omega}  \big\{ f(x) \doteq \frac{1}{m}\sum_{i=1}^{m} f^i (x)  \big\},
\end{equation}
where $f^i :\mathbb{R}^d \to \mathbb{R}$ for all $i$. This situation arises quite frequently in machine learning in the context of stochastic optimization While it may be desirable to minimize (\ref{prob1}), such a goal is sometimes untenable when one does not have access to the law $P$ or when one cannot draw from an infinite population sample set. In that case,  a practical approach is to instead seek the solution of a problem that involves an estimate of the the expectation in (\ref{prob1}) giving rise to (\ref{prob}). Then, $f^i(x) = f(x,\xi_i)$ in (\ref{prob}) for some realization $\xi_i$ of $\xi$.

Optimization problems of the form (\ref{prob}) usually have a large $m$, which makes finding a possible solution by first order methods a computationally intensive and time consuming task. A possible remedy is to use a stochastic approach.  Also, because of high dimensions, they are generally beyond the capability of second order methods due to extremely high iteration complexity (see \cite{bot}). Note that (\ref{prob}) is a constrained optimization problem, so the conventional approaches (e.g., projected gradient descent) to solve it require a projection on the constraint set at every step of the algorithm. This can be quite an expensive operation and may make the problem intractable.

Recently there has been a lot of interest in the Frank-Wolfe method (FW) \cite{frank}, also known as conditional gradient, for solving constrained optimization problems. In many problems where the domain $\Omega$ admits fast linear optimization while having slow projections, FW has the potential to be much more efficient than projection based algorithms. Such problems include linear regression with constraints, multiclass classification, matrix completion, etc.\\

\textbf{Contributions : } 
Our contributions are as follows. We suggest three algorithms related to Frank-Wolfe method to deal with (\ref{prob1}) and (\ref{prob}) : 
\begin{itemize}

\item The first algorithm is based on the dual averaging technique of \cite{nest}, \cite{xiao} and deals with (\ref{prob1}). The essential idea of \cite{nest}  is to use a weighted average of all the gradients till the current step  instead of just the gradient at the current iterate. This is done to overcome a significant disadvantage of the latter approach. Since a diminishing time step is used in most gradient related algorithms, this implies new gradient computations enter the algorithm with decreasing weights (see (1.6), \cite{nest}). We leverage this idea to propose a stochastic version of FW, where $\mathcal{O}(\frac{1}{\epsilon^3})$ gradient evaluations are required while having a convergence rate of $\mathcal{O}(\frac{1}{\epsilon})$.\footnote{By this we mean $\mathcal{O}(\frac{1}{\epsilon})$ steps are needed to reach any $\epsilon$-accurate solution, so that $f(x) - \min_{y \in \Omega} f(y) \leq \epsilon$ for any output $x$ of the algorithm.}   

\item We also propose an incremental gradient approach to solve (\ref{prob}). This requires only  $\mathcal{O}(\frac{1}{\epsilon})$ gradient evaluations while maintaining the same convergence rate. We remark here that under our assumptions (convex and smooth), this is the best known convergence rate with the minimal number of gradient evaluations to deal with a sum of functions form. The incremental method involves storing the last $m$ gradient evaluations and using a running average of these $m$ evaluations instead of just the gradient at the current iterate.  The essential idea why the incremental version works is that because of a small time step which is diminishing, the last $m$  iterates generated by the algorithm are not very different from each other. Together with the smoothness of the function, this implies a good enough approximation of the running average to the gradient at the current iterate. We also remark here the proof of the main convergence result (Theorem 3) is remarkably easy and may have application elsewhere.

\item Finally, we also propose a distributed version of FW with dual averaging. The proposed algorithm is inspired from the ideas of \cite{duc}, where distributed dual averaging is developed for the proximal gradient method. This algorithm also achieves optimal convergence rates which could be achieved for a distributed implementation of FW.
\end{itemize}

\textbf{Related Work : }  The classical Frank-Wolfe method using a line search was originally proposed in  \cite{frank}. The convergence rate for a smooth convex function with a polyhedral domain was shown to be $\mathcal{O}(\frac{1}{\epsilon})$. This rate was later improved  in \cite{wolf} for such domains  by suggesting a variant know as the "away step FW". A more recent treatment of the method was given in  \cite{jag}, where it was extended to atomic domains. For the classical Frank-Wolfe algorithm, \cite{beck} showed a linear rate for the special case of quadratic objectives when the optimum is in the strict interior of the domain, a result also proved earlier in  \cite{gue}. Building upon the work of \cite{levit}, it was shown in \cite{garb} that for strongly convex sets, FW achieves a rate of $\mathcal{O}(\frac{1}{\sqrt{\epsilon}})$. We remark here that any improvement in the convergence rate of FW depends upon the geometric properties (being a polytope, strong convexity etc.) of the constraint set. A primal-dual averaging algorithm for deterministic optimization was studied in \cite{lan2}. We also mention \cite{aber} where a deterministic algorithm similar to dual averaging was presented in a game theoretic setting of online learning. We also mention the related the work \cite{lu} in the context of stochastic FW, where the authors consider a stochastic substitute version of FW in the context of structured convex optimization.

Compared to the deterministic version, the stochastic version is relatively less studied. One of the main works to study FW in an online learning setting is \cite{haz}. The vanilla version of stochastic FW requires $\mathcal{O} (\frac{1}{\epsilon^3})$ stochastic gradient evaluations to achieve an $\epsilon$-accurate optimum. The convergence rate is $\mathcal{O}(\frac{1}{\epsilon})$. The alternative to the vanilla version suggested in \cite{haz} is based on combing gradient sliding techniques \cite{lan} with variance reduction algorithms \cite{john}, \cite{mah}. Although the convergence rate remains the same, the number of gradient evaluations is brought down significantly. In \cite{red}, a similar approach (variance reduction) is taken to handle non-convex functions. The convergence rate and gradient evaluations required are both of the order $\mathcal{O} (\frac{1}{\epsilon^2})$. 

The literature on distributed algorithms is vast, mostly building upon the seminal work of \cite{tsit1}.  The most relevant works for our purposes are \cite{bell} and \cite{wai}. In \cite{bell}, a distributed FW is proposed in order to solve the associated optimization problem in the  setting where the elements to be combined are not centrally located but spread over a network. A decentralized version (where at each iteration an approximate average is obtained ) is studied in \cite{wai}. Furthermore, a  much more communication efficient version is also proposed which utilizes the gradient tracking concept of \cite{ned}. The distributed version we study here is inspired by \cite{duc}, which studies distributed dual averaging for proximal gradient methods. Broadly speaking, this work also leverages two time scale algorithm theory whose braod range of applications can be found in \cite{s2,s3,s4,s5,s6,s9,s10} maong others.

\section{Background }

In this section we discuss the basics of Frank-Wolfe, dual averaging and distributed algorithms \cite{s1,s7,s8}. 

\subsection{Frank-Wolfe }
We first review the standard Frank-Wolfe algorithm, the pseudo code for which is given in Algorithm 1.

\begin{algorithm}[H]
\textbf{Input :} Objective Function  $f$ ; parameter $\alpha_k$ \\

\textbf{Initialize :} $x_0 \in \Omega$, $w_0 = \text{argmin}_{x \in \Omega} \nabla f(x_0)^T x$. \\

\textbf{For} $k =1,2...$ \textbf{do} : \\ 
\begin{enumerate}

\item Compute  the gradient $\nabla f(x_{k})$. Call the linear optimization oracle to compute :
$$w_k \in  \text{argmin}_{x \in \Omega} \, \langle \nabla f(x_{k}) , x \rangle. $$

\item Set $$ x_{k+1} = (1-\alpha_k)x_{k} + \alpha_k w_{k}.$$ 

\end{enumerate}

\textbf{end for }\\

\caption{Classical Frank-Wolfe (FW)}
\end{algorithm}

At each instant $k$, the algorithm makes a call to a linear optimization (LO) oracle, followed by a convex averaging step of the current iterate with the output of the LO oracle (Step b). Note that as long as $x_0$ is initialized in the interior of $\Omega$, the sequence $\{ x_k\}_{k \geq 1 }$ is guaranteed to stay inside $\Omega$.  Since the algorithm substitutes the projection problem (which can be thought of as a constrained quadratic minimization problem) with a linear optimization one, FW is sometimes referred to as a "projection-free" algorithm. This can make it very desirable for structured constraints ( see \cite{lac}). The deterministic Frank-Wolfe has a convergence rate of $\mathcal{O} (\frac{1}{k})$ (Theorem  1, \cite{jag}). Besides the advantages of a low iteration cost and ease of implementation, the iterates generated by FW enjoy many structural properties. In particular, since the iterates $\{x_k\}$ can be written as a convex combination of a smaller number of extreme points of $\Omega$, sparsity and low rank (for matrix constraints) are preserved at every step (Section 3, \cite{jag}).

\subsection{Dual Averaging}  We review the idea of  dual averaging in the context of mirror descent. Mirror descent is based on a proximal function $\psi : \Omega \to \mathbb{R}$ assumed to be 1-strongly convex with respect to the Euclidean norm $\|\cdot\|$ :
$$
\psi(y) \geq \psi(x) + \langle  \nabla \psi (x),y-x\rangle + \frac{1}{2}\|x-y\|^2.
$$
An example of such a proximal function is the quadratic function $\psi(x)=\frac{1}{2} \|x\|^2$. Define the mapping :
$$
\Pi^{\Psi}_{\Omega}(g,\alpha) = \text{arg min}_{ x  \in \Omega} \Big\{ \langle g,x\rangle + \frac{1}{\alpha} \psi(x)\Big\}. 
$$
We state the  mirror descent algorithm with  dual averaging in Algorithm 2.

\begin{algorithm}[H]
\textbf{Input :} Objective Function  $f$ ; parameter $\alpha_k$ \\

\textbf{Initialize :} Set $x_0 = \text{argmin}_{x \in \Omega}\, \psi(x)$ and $g_0 = 0 $ .\\

\textbf{For} $k =1,2...$ \textbf{do} : At each node $i$, \\ 

\begin{enumerate}

\item Compute  the gradient $\nabla f(x_{k})$.\\ 

\item Update the average gradient :
$$
g_{k} = \frac{k-1}{k} g_{k-1} + \frac{1}{k} \nabla f(x_{k}).
$$

\item Set $$ x_{k+1} = \Pi^{\Psi}_{\Omega}(g_{k},\alpha_k). $$ 

\end{enumerate}

\textbf{end for }

\caption{ Mirror descent with Dual Averaging (FW)}
\end{algorithm}
 We note that,
$$
g_k =\frac{1}{k}\sum_{i=1}^{k} \nabla f(x_i).
$$ 
This implies that we use the average gradient upto time $k$ in Step (b) of Algorithm 2. For the stochastic setting, instead of using the entire sum present in $f(\cdot)$,\footnote{Recall that $f(\cdot) = \sum_{i=1} ^m f^i(\cdot)$.} one can randomly sample $f^i$ and use $\nabla f^i(\cdot)$ instead of $\nabla f(\cdot)$ in Step (c). We refer the reader to \cite{xiao} for more details on dual averaging in an online setting. We remark that the above algorithm is equivalent to the standard projected gradient descent with the  average gradient when $\psi(x)=\frac{1}{2} \|x\|^2 $ (Proposition 2.1, \cite{ber}).

\subsection{Distributed Algorithms } Suppose we have a network of $m$ agents indexed by $1, ..., m$. We associate with each agent $i$, the function $f^i : \mathbb{R}^d \to \mathbb{R}$ and a global convex constraint set $\Omega$.\footnote{In many applications, one may also have $f^i(\cdot) = f(\cdot)\,\forall i$.} The (global) function which we aim to minimize is $f:\mathbb{R}^d \to \mathbb{R}$.

 Let the communication network be modelled by a static undirected graph
$\mathcal{G=}\{\mathcal{V},\mathcal{E}\}$ where $\mathcal{V}=\{1,...,m\}$
is the node set and $\mathcal{E\subset\mathcal{V}\mathcal{\times}\mathcal{V}}$
is the set of links $(i,j)$ indicating that agent $j$ can send information
to agent $i$. All of the arguments presented here can be extended
to a time-varying graph under suitable assumptions. Here we deal only with
a static network for ease of notation.

We associate with the network a non-negative weight matrix
 $Q = [[q_{ij}]]_{i,j \in \V}$ such that
 $$q_{ij} > 0 \Longleftrightarrow (i,j)\in\mathcal{E}.$$
\textbf{Assumption 1:}\\

i) {[}\textit{Double Stochasticity}{]} $\mathbf{1}^{T}Q=\mathbf{1}^{T}$
and $Q\mathbf{1}=\mathbf{1}$.\\

ii) {[}\textit{Irreducibility and aperiodicity}{]} We assume that the underlying graph is irreducible, i.e., there is a directed path from any  node to any other node, and  aperiodic, i.e., the g.c.d.\ of lengths of all paths from a node to itself is one. It is known that the choice of node in this definition is immaterial. This property can be guaranteed, e.g.,  by making $q_{ii}>0$ for some $i$.\\

The objective of distributed optimization is to minimize $f(\cdot)$ subject to staying in the constraint $\Omega$  while simultaneously achieving consensus, i.e.
\begin{align*}
\text{minimize} & \qquad \sum_{i=1}^{m} f^i(x^i)\\
\text{subject to} & \qquad \,x^i \in \Omega\\
& \qquad \sum_{j \in \mathcal{N}_i}q_{ij} \|x^i -x^j\|^2 =0\,\,\,\forall i.
\end{align*}
The most popular way to solve this is studied in for the unconstrained case in \cite{tsit1}   and subsequently in \cite{ned2010} for the constrained case. It involves the following two steps :

\begin{enumerate}
\item[(S1)] {[}\textit{Consensus Step}{]} This step involves local averaging at each node and is aimed at achieving consensus, 
$$v^i_k = \sum_{j\in\N(i)} q_{ij}x^j_k. $$

\item[(S2)] {[}\textit{Gradient Descent Step}{]} This step is the gradient descent part aimed at minimizing $f^i$ at each node :
$$x^i_{k+1} = P_{\Omega}(v^i_k - a_k \nabla f^i(v^i_k) ), $$
where $P_{\Omega}(\cdot)$ is the projection on the set $\Omega$ and $a_k$ is a positive scalar decaying at a suitable rate to zero. 

\end{enumerate}
The convergence properties of the above algorithm have been extremely well studied for convex as well as the non-convex case.

\section{The algorithms }
We assume here that the function $f$ in (\ref{prob1}) (and each $f^i$ in (\ref{prob})) is convex and $L$-smooth in $\mathbb{R}^d$. By smoothness, we mean that the gradient is $L$-Lipschitz continuous, i.e.
$$
\| \nabla f(x) -\nabla f(y)\| \leq L \|x-y\|.
$$
For a general $f : \mathcal{X} \to \mathbb{R} $, the norm on the LHS will be the dual norm. Together with the convexity, the $L$-smoothness of $f$ implies :
\begin{equation} \label{smooth}
f(y)  \leq f(x) + \langle \nabla f(x), y-x  \rangle + \frac{L}{2} \| x-y\|^2 .
\end{equation}
We use $L_f(x,\cdot)$ to denote the linear approximation of $f(\cdot)$ at the point $x$ :
\begin{equation}\label{linear}
L_f(x,y) :=f(x) + \langle \nabla f(x), y-x  \rangle.
\end{equation}
Let $C_{\Omega}$ denote the following bound :
\begin{equation}\label{bound}
C_{\Omega} := \sup_{x,y \in \Omega} \| x- y\|.
\end{equation}
\subsection{Frank-Wolfe with stochastic dual averaging (FW-SDA)}
Our first algorithm combines the Frank-Wolfe algorithm with stochastic dual averaging  to solve (\ref{prob1}). We assume here that we have access to a stochastic first order oracle (SFO). The oracle takes a point $x$ and returns an unbiased sample $\nabla f(x,\xi')$, where $\xi'$ is a sample drawn i.i.d from $P$. By an unbiased sample $\nabla f(x,\xi')$, we mean that 
$$
\mathbb{E} [\nabla f(x,\xi')] = \nabla f(x).
$$
The essential idea here is to generate enough of these random samples at each iteration and  use it to update the average gradient employed in a standard dual averaging scheme. The pseudo code is provided in Algorithm 3. In addition to the  sequence $\{x_k\}$ and $\{w_k\}$ which are present in the Frank-Wolfe algorithm, Algorithm 3 has two additional sequences $\{z_k\}$ and $\{g_k\}$. The auxiliary sequence $\{z_k\}$  is again standard for Nesterov's algorithm. The sequence $\{g_k\}$ keeps track of the average gradient upto time $k$ and constitutes the dual averaging part of the algorithm. We note that $g_k$ can be written as :

\begin{equation}
g_k = \frac{1}{\sum_{i=1}^{k}\beta_i} \sum_{i=1}^{k} \beta_i \nabla_ {k} (z_{i-1}), 
\end{equation}
where $\nabla_k$ is as in (\ref{unbiased}).

\begin{algorithm}[H]
\textbf{Input :} Initial point $x_0$ ; parameters $\alpha_k$, $p_k$ and $\beta_k$.\\

\textbf{Initialize :} $x_0 \in \Omega$, $w_0 = \text{argmin}_{x \in \Omega} \langle \nabla f(x_0) , \,x\rangle$ and $g_0 =0$.\\

\textbf{For} $k =1,2...$ \textbf{do} :\\ 
\begin{enumerate}

\item  Set
 \begin{equation}\label{SFW-1}
z_{k-1} = (1-\alpha_k)x_{k-1} + \alpha_k w_{k-1}.
\end{equation}

\item  Draw $p_k$ i.i.d samples ${\xi_1,....,\xi_{p_k}}$ according to the distribution $P$ independent of the past history $\mathcal{F}_{k-1}$ and set 
\begin{equation}\label{unbiased}
\nabla_k (z_{k-1}):= \frac{1}{p_k}  \sum_{j=1}^{p_k} \nabla f(z_{k-1}, \xi_j)
\end{equation}

\item Update the average weighted gradient  :
\begin{equation}\label{SFW-2}
g_k = \frac{1}{ \sum_{i=1}^{k}\beta_i} \Big\{ \big( \sum_{i=1}^{k-1} \beta_i \big)g_{k-1} + \beta_k  \nabla_k (z_{k-1}) \Big\}.
\end{equation}

\item Call the Linear Optimization (LO) oracle to compute :
$$w_k \in  \text{argmin}_{x \in \Omega} \, \langle g_k , x \rangle. $$

\item Set 
\begin{equation} \label{SFW-3}
 x_{k} = (1-\alpha_k)x_{k-1} + \alpha_k w_{k}.
\end{equation}
\end{enumerate}

\textbf{end for }\\

\caption{Stochastic Frank-Wolfe with Dual averaging (FW-SDA)}
\end{algorithm}

Let $L_{f^{k}}$ denote the following (with $\nabla_k(x)$ as in (\ref{unbiased})),
\begin{equation}\label{f^k}
L_{f^{k}} (x,y ):= f(x) + \langle \nabla_k(x), y-x  \rangle.
\end{equation}
for any $x, y \in \Omega $. We have, 
\begin{equation*}
\mathbb{E} [L_{f} (x,y )-L_{f^k}(x,y) ] = \mathbb{E} [ \langle \nabla f(x)- \nabla_k(x), y-x  \rangle ].
\end{equation*}
From Cauchy Schwarz inequality,  
\begin{equation}\label{ftof^kbound}
\mathbb{E} [L_{f} (x,y )-L_{f^k}(x,y) ] = C_{\Omega} \mathbb{E} [\|  \nabla f(x)-\nabla_k (x)\|  ].
\end{equation}
where  $C_{\Omega} = \sup_{x,y \in \Omega } \| x-y \|$.
Step (d) of Algorithm 3  implies 
\begin{equation}\label{min}
w_k  \in \text{argmin}_w\,\, \Phi_k (w), 
\end{equation}
where  $\Phi_k (w)$ is the accumulated linear model till time $k$, i.e., $\Phi_k (w) =0 $ if $k=0$ and 
\begin{equation}\label{phi}
\Phi_k (w) : = \frac{1}{\sum_{i=1}^{k}\beta_i} \sum_{i=1}^{k} \beta_i L_{f^{i}} (z_{i-1} ; w ), k \geq 1 
\end{equation} 
with $L_{f^{k}}(\cdot,\cdot)$ defined in (\ref{f^k}). The next theorem gives the rate of convergence of Algorithm 3 and builds upon the proof of its deterministic counterpart \cite{lan2}.

\begin{thm} \label{DAthrm}
Let $\alpha_k = \frac{2}{k+1} $. If the parameter $\beta_k$ is chosen such that
\begin{equation}\label{cond}
\alpha_k  = \frac{\beta_k}{\sum_{i=1}^{k}\beta_i},
\end{equation} 
 and the number of samples in $p_k$ in (\ref{unbiased}) is such that the following bound holds,  
\begin{equation}\label{variance-bound}
\mathbb{E}[\| \nabla_k(x) - \nabla f(x)\| ]\leq \frac{1}{k}, 
\end{equation}
for any $x\in \mathbb{R}^d$, then Algorithm 3 ensures that
$$ \mathbb{E}[f(x_k) -f(x^*)]   = \mathcal{O}(\frac{1}{k}) $$
 for any $k>0$. \\

\end{thm}

\begin{proof}
(i) We have from  (\ref{smooth}),
\begin{align*}
f(x_k) &\leq L_{f} (z_{k-1}; x_k ) + \frac{L}{2} \|x_k -z_{k-1} \| ^2,\\ 
&= (1 - \alpha_k) L_{f} (z_{k-1}; x_{k-1} ) + \alpha_k L_{f}  (z_{k-1}; w_{k} )+ \frac{L}{2} \|x_k -z_{k-1} \| ^2,
\end{align*} 
where we have used (\ref{SFW-3}) and the linearity of $L_{f}(\cdot;\cdot)$ in the second argument.  Subtracting (\ref{SFW-1}) from (\ref{SFW-3})  we have,
$$\| x_k -z_{k-1} \| = \alpha_k \| w_k - w_{k-1}\| .$$
Since $C_{\Omega} = \sup_{x,y \in \Omega } \| x-y \|$, we get 
\begin{equation*}
f(x_k)\leq (1 - \alpha_k) L_{f} (z_{k-1}; x_{k-1} ) + \alpha_k L_{f}  (z_{k-1}; w_{k} ) + \frac{L}{2}\alpha^2_k  C_{\Omega}^2.
\end{equation*} 
Using the convexity of $f(\cdot)$ in the first term in the RHS in the above we get,
\begin{equation*}\label{eq*}
f(x_k) \leq (1 - \alpha_k) f(x_{k-1})  + \alpha_k L_{f}  (z_{k-1}; w_{k} ) + \frac{L}{2}\alpha^2_k  C_{\Omega}^2 .
\end{equation*}
Taking expectation we have,
\begin{equation}\label{eq-22}
 \mathbb{E}[ f(x_k) ]  \leq (1 - \alpha_k) \mathbb{E}[f(x_{k-1})]  +  \alpha_k \mathbb{E}[ L_{f^k}  (z_{k-1}; w_{k} ] + \alpha_k \mathbb{E} [\{ L_{f}  (z_{k-1}; w_{k} )- L_{f^k}  (z_{k-1}; w_{k} )\}] + \frac{L}{2}\alpha^2_k  C_{\Omega}^2 .
\end{equation}
We note that the third term in the RHS of the above inequality can be bounded using (\ref{ftof^kbound}) and (\ref{variance-bound}) as 
$$
\mathbb{E} [\{ L_{f}  (z_{k-1}; w_{k} )- L_{f^k}  (z_{k-1}; w_{k} )\}]\leq \frac{C_{\Omega}}{k} \leq C_{\Omega}\alpha_k, 
$$
so that 
\begin{equation}\label{eq*}
 \mathbb{E}[ f(x_k) ]  \leq (1 - \alpha_k) \mathbb{E}[f(x_{k-1})]  +  \alpha_k \mathbb{E}[ L_{f^k}  (z_{k-1}; w_{k} ] + C_{\Omega} \alpha_k^2 + \frac{L}{2}\alpha^2_k  C_{\Omega}^2 .
\end{equation}
Let $ \bar{\beta}_k = \sum_{i=1}^{k}\beta_i$. We next bound the second term in the RHS of (\ref{eq-22}). We have from (\ref{phi}),
 $$  \bar{\beta}_{k} \Phi_k (w_k)   = \beta_k L_{f^{k}} (z_{k-1} ; w_k ) + \bar{\beta}_{k-1} \Phi_{k-1} (w_k) . $$
From (\ref{min}), $w_{k-1}  \in \text{argmin}_w\,\, \Phi_{k-1} (w) $, so that $ \Phi_{k-1}(w_{k-1}) \leq   \Phi_{k-1} (w_k)$. This gives 
\begin{align}\label{aux-1}
\bar{\beta}_{k}  \Phi_k (w_k)  \geq \beta_k L_{f^{k}} (z_{k-1} ; w_k ) + \bar{\beta}_{k-1}  \Phi_{k-1} (w_{k-1}) .
 \end{align} 
We also have from (\ref{cond}) that $\alpha_k \beta_k^{-1}=\bar{\beta}_k^{-1}$. So,
\begin{align*}
\alpha_k \beta_{k}^{-1} [\bar{\beta}_{k}  \Phi_k (w_k) - \bar{\beta}_{k-1}  \Phi_{k-1} (w_k)] &=  \Phi_k (w_k) -\alpha_k \beta_{k}^{-1}  \bar{\beta}_{k-1}  \Phi_{k-1} (w_{k-1}),\\
&=  \Phi_k (w_k) - \bar{\beta}_{k-1}  \bar{\beta}^{-1}_{k}\Phi_{k-1} (w_{k-1}).
\end{align*} 
Noting that $\bar{ \beta}_{k-1} = \bar{ \beta}_{k} - \beta_k  $, we have $ \bar{\beta}_{k-1}  \bar{\beta}^{-1}_{k}=( 1-  \bar{ \beta}^{-1}_{k} \beta_k)=(1-\alpha_k)$ (using (\ref{cond})). So, 
\begin{align*}
\alpha_k \beta_{k}^{-1} [\bar{\beta}_{k}  \Phi_k (w_k) -\bar{\beta}_{k-1}  \Phi_{k-1} (w_k)]  &=  \Phi_k (w_k) -  (1-\alpha_k ) \Phi_{k-1} (w_{k-1}).
\end{align*} 
Using the above equality in  (\ref{aux-1}), we have 
\begin{equation*}
 \alpha_k \mathbb{E}[ L_{f^{k}}(z_{k-1},w_k) ] \leq   \mathbb{E}[ \Phi_k (w_k)  -  (1-\alpha_k ) \Phi_{k-1} (w_{k-1})].
\end{equation*}
Then using the above in (\ref{eq*}),
\begin{equation*}
  \mathbb{E}[ f(x_k) ]  \leq (1 - \alpha_k) \mathbb{E}[f(x_{k-1})]+   \mathbb{E}\big[ \Phi_k (w_k) -  (1-\alpha_k )\Phi_k (w_{k-1})\big]+  C_{\Omega} \alpha_k^2 + \frac{L}{2}\alpha^2_k  C_{\Omega}^2.
\end{equation*}
Rearranging, we get
\begin{equation}\label{final}
\mathbb{E}[ f(x_k)-\Phi_k (w_k)   ] \leq (1 - \alpha_k) \mathbb{E}[(f(x_{k-1}) - \Phi_k(w_{k-1})] +  C_{\Omega} \alpha_k^2+ \frac{L}{2}\alpha^2_k  C_{\Omega}^2 .
\end{equation}
Set $ \nu_k =  \mathbb{E}[ f(x_k)- \Phi_k (w_k)   ]  $ and  .define 
\begin{equation*}
   \Gamma_k :=
    \begin{cases}
      (1-\alpha_k)  \Gamma_{k-1}, & \text{if}\ k\geq 2 \\
      1, &  \text{if}\ k=1.
    \end{cases}
\end{equation*}
With the above notation, we can write (\ref{final})  as :
\begin{align*}
\nu_k &\leq (1 - \alpha_k) \nu_{k-1} +  \frac{L}{2}\alpha^2_k  C_{\Omega}^2 + C_{\Omega} \alpha_k^2 ,\\
\Longrightarrow \frac{\nu_k }{\Gamma_k} &\leq  \frac{ \nu_{k-1}}{\Gamma_{k-1}} + \frac{L}{2}\frac{\alpha^2_k}{\Gamma_k}  C_{\Omega}^2 +  C_{\Omega}\frac{\alpha^2_k}{\Gamma_k} .  \,\,\,\,\, 
\end{align*}
A simple calculation gives :
$$
\Gamma_k = \frac{2}{k(k+1)} \,\,\,\text{       and        }\,\,\,\,  \frac{\alpha_k^2}{\Gamma_k}\leq 2.
$$
To conclude the proof we note that,
\begin{align*}
\nu_k &\leq \Gamma_k \Big\{ \frac{\nu_1}{\Gamma_1} + \sum_{i=1}^{k} \big( \frac{L}{2}  C_{\Omega}^2 + C_{\Omega} \big)  \Big( \frac{\alpha^2_i}{\Gamma_i} \Big) \Big\},\\
\nu_k &\leq \frac{2}{k(k+1)}\Big\{ \frac{\nu_1}{\Gamma_1} + k(L C_{\Omega}^2 + 2C_{\Omega})\Big\},\\
& = \mathcal{O}(\frac{1}{k}).
\end{align*}
From (\ref{min}) we have,
$$
\Phi_k (w_k)  \leq \Phi_k (w)\,\,\forall w\in\Omega. 
$$
Then, by the definition of $\Phi_k (w_k) $,
$$
\mathbb{E} [\Phi_k (w_k)]  \leq  \frac{1}{\sum_{i=1}^{k}\beta_i}\mathbb{E} [ \sum_{i=1}^k \beta_i  f(z_{i-1}) + \beta_i  \langle \nabla_i (z_{i-1}) , w - z_{i-1} \rangle  ].  
$$
From the unbiasedness of $\nabla f(x,\xi)$ and the independence of $z_{i-1}$ and $\xi$ (see Step (b)), we have 
$$
\mathbb{E} [\Phi_k (w_k)]   \leq   \frac{1}{\sum_{i=1}^{k}\beta_i}  \sum_{i=1}^k \{ \beta_if(z_{i-1}) +   \beta_i\langle \nabla f(z_{i-1}) , w - z_{i-1} \rangle  \},
$$
which gives from the convexity of $f$,
$$
\mathbb{E} [\Phi_k (w_k)] \leq f(w).
$$
We have for $w=x^*$ in the above, 
$$
\mathbb{E}[ \Phi_k(w_k) ]   \leq    f(x^*) ,
$$
which gives
$$
 \mathbb{E}[ f(x_k)-f (x^*)   ]   \leq  \mathbb{E}[ f(x_k)- \Phi_k (w_k)   ] =  \nu_k = \mathcal{O}\Big(\frac{1}{k}\Big).
$$
\end{proof}
\begin{rem}
The condition (\ref{variance-bound}) can be satisfied by setting $p_k = k^2$ and the condition (\ref{cond}) can be satisfied by using $\beta_k = k,\,k\geq1$.
\end{rem}

\subsection{ Incremental Method :} 
In this section we consider an incremental approach to Frank-Wolfe for solving (\ref{prob}). For a survey of incremental methods in the context of gradient descent and proximal algorithms, we refer the reader to \cite{ber}. The driving idea behind incremental methods for problems of form (\ref{prob}) is to keep track of the running average of the last $m$-gradient evaluations. Let $[k]_m: = (k \text{ modulo }m) +1$. These $m$ gradient evaluations can be constructed in a number of ways. The most obvious way is to go through each $f^i$ in a cyclic manner, so that at each instant $k$, we have the following average :
$$
d_k = \sum_{i=1}^{m} \nabla f^{[k]_m} (x_{k-i+1}).
$$
Another possible technique, popular in neural network training practice, is to reshuffle
randomly the order of the component functions after a cycle of $m$ steps. Algorithm 4 gives the pseudo-code for incremental Frank-Wolfe. We take the cyclic approach here. The convergence analysis for the random shuffling algorithm differs only in trivial ways from the cyclic one.

\begin{algorithm}[H]
\textbf{Input :} Objective Function  $f=\frac{1}{m}\sum_{i=1}^m f^i$ ; parameters $\alpha_k$ \\

\textbf{Initialize :} $x_0 \in \Omega$, $w_0 = \text{argmin}_{x \in \Omega} \langle \nabla f(x_0),\,x\rangle$ and $d_0 =0$.\\

\textbf{For} $k =0,1,...$ \textbf{do} :\\ 
\begin{enumerate}

\item  Compute the gradient $\nabla f^{[k]_{m}}(x_{k})$.\\

\item Update the aggregated gradient :
$$d_k = \frac{1}{m}  \Big\{  md_{k-1} - \nabla f^{[k]_{m}}(x_{k-m}) + \nabla f^{[k]_{m}}(x_{k}) \Big\} $$.

\item Call the LO oracle to compute :

$$w_k \in  \text{argmin}_{x \in \Omega} \langle d_k ,  x \rangle $$.

\item Set \begin{equation}\label{stepd} x_{k+1} = (1-\alpha_k)x_{k} + \alpha_k w_{k}.
\end{equation}
\end{enumerate}

\textbf{end for }\\

\caption{Incremental Frank-Wolfe }
\end{algorithm}

\begin{thm} \label{DAthrm}

With $\alpha_k = \frac{2}{k+1} $, Algorithm 4 ensures that 

$$ f(x_k) -f(x^*)  = \mathcal{O}\Big(\frac{1}{k}\Big)   $$
 
for any $k>0$. 
\end{thm}

\begin{proof}
We have 
\begin{align*}
f(x_{k+1}) &\leq f(x_{k}) + \langle  \nabla f(x_{k}), x_{k+1} - x_{k} \rangle  +\frac{L}{2} \|x_{k+1} -x_{k} \| ^2 \,\,\,(\text{from } (\ref{smooth}))\\
&=  f(x_{k}) +  \alpha_k \langle \nabla f(x_{k}), w_k - x_{k} \rangle  +  \frac{L}{2}\alpha_k^2 \|w_{k} -x_{k} \| ^2 \,\,\, (\text{from } (\ref{stepd}))\\
&=  f(x_{k}) + \alpha_k \langle  \nabla f(x_{k}) - d_k , w_k - x_{k} \rangle  +   \alpha_k  \langle  d_k, w_k - x_{k} \rangle + \frac{L}{2}\alpha^2 \|w_{k} -x_{k} \| ^2.
\end{align*} 
From the optimality of $w_k$, we have 
\begin{equation}\label{inc-0}
f(x_{k+1})  \leq   f(x_{k}) + \alpha_k \langle  \nabla f(x_{k}) - d_k , w_k - x_{k} \rangle  +
 \alpha_k  \langle  d_k, x^* - x_{k} \rangle + \frac{L}{2}\alpha^2 C_{\Omega}^2.
\end{equation}
By the definition of $d_k$,
\begin{align*}
\| d_k - \nabla f(x_{k})\| &= \frac{1}{m} \| \sum_{i=k-m+1}^{k} \nabla  f^{[i]_m} (x_i) - \sum_{i=1}^{m} \nabla f ^i (x_{k})\|.
\end{align*} 
Since $\sum_{i=1}^{m} \nabla f^i (x_{k}) = \sum_{i=k-m+1}^{k} \nabla  f^{[i]_m} (x_{k}) $. we have 
\begin{align*}
\| d_k - \nabla f(x_{k})\| &= \frac{1}{m} \| \sum_{i=k-m+1}^{k} \big( \nabla  f^{[i]_m} (x_i) - \nabla  f^{[i]_m}(x_{k}) \big)\|.
\end{align*}
From the smoothness of $f^i(\cdot)$ for all $i$, we have
\begin{equation}\label{inc-1}
\| d_k - \nabla f(x_{k})\| \leq \frac{L}{m}\sum_{i=k-m+1}^{k} \|  x_i- x_{k} \| 
= \frac{L}{m} \sum_{i=1}^{m} \|x_{k-i+1} - x_{k} \|. 
\end{equation}
We note from (\ref{stepd}), that for $1  \leq i \leq  m-1$,
\begin{equation}
x_{k} - x_{k-i} = \sum_{j=1}^{i} \alpha_{k-j} (w_{k-j} - x_{k-j} )
\end{equation}
so that
\begin{equation}\label{inc-2}
\| x_{k} -x_{k-i} \| \leq \sum_{j=1}^{i} \alpha_{k-j} \|w_{k-j} - x_{k-j} \| \leq C_{\Omega} \sum_{j=1}^{i} \alpha_{k-j}.
\end{equation}
For $k>2m$, we have  for $1\leq j \leq m$ $$\frac{\alpha_{k-j}}{\alpha_k}\leq \frac{\alpha_{k-m}}{\alpha_k} \leq 2,$$ so that
$$
\sum_{j=1}^{i} \alpha_{k-j} \leq  \sum_{j=1}^{i} \Big( \frac{\alpha_{k-j}}{\alpha_k}\Big) \alpha_k \leq 2\sum_{j=1}^{i}  \alpha_k \leq 2i \alpha_k. 
$$
 We use the above bound in (\ref{inc-2}) to get
$$
\| x_{k} -x_{k-i} \| \leq 2C_{\Omega} \times  i \alpha_k .
$$
Using this in (\ref{inc-1}),
\begin{equation}\label{111new}
\| d_k - \nabla f(x_{k-1})\| \leq  \frac{2LC_{\Omega}  \alpha_k}{m} \sum_{j=1}^{m} j= LC_{\Omega}\alpha_k(m+1).
\end{equation}
We use the above bound in (\ref{inc-0}) :
\begin{align*}
f(x_{k+1})  &\leq   f(x_{k}) + \alpha_k  \langle  d_k, x^* - x_{k} \rangle   + \alpha_k \langle  \nabla f(x_{k}) - d_k , w_k - x_{k} \rangle + \frac{L}{2}\alpha_k^2 C_{\Omega}^2,\\
& \leq  f(x_{k}) + \alpha_k \langle  \nabla f(x_{k}) , x^* - x_{k} \rangle  +  \alpha_k  \langle   \nabla f(x_{k})- d_k, w_k - x^* \rangle + \frac{L}{2}\alpha_k^2 C_{\Omega}^2,\\
& \leq  f(x_{k}) + \alpha_k \langle  \nabla f(x_{k}) , x^* - x_{k} \rangle  +  \alpha_k C_\Omega  \|  \nabla f(x_{k})- d_k \|+ \frac{L}{2}\alpha_k^2 C_{\Omega}^2,\\
& \leq f(x_{k}) + \alpha_k \langle  \nabla f(x_{k}) , x^* - x_{k} \rangle   + \frac{LC_\Omega^2\alpha^2_k(m+1)}{2}  + \frac{L}{2}\alpha_k^2 C_{\Omega}^2.
\end{align*}
This gives from the convexity of $f(\cdot)$,
\begin{equation*}
f(x_{k+1}) \leq f(x_{k}) + \alpha_k \Big(  f(x^*)-f(x_{k})\Big)  
+\frac{LC_\Omega^2\alpha^2_k}{2} (m+2),
\end{equation*}
so that 
\begin{equation*}
f(x_{k+1}) -f(x^*) \leq (1-\alpha_k)(f(x_{k})-f(x^*)) 
+\frac{LC_\Omega^2\alpha^2_k}{2} (m+2).
\end{equation*}
Setting $\nu_k =f(x_k) -f(x^*)  $ we have,
\begin{equation*}
\nu_{k+1}  \leq (1-\alpha_k)\nu_k +    \mathcal{O}(\alpha_k^2).
\end{equation*}
Proceeding  the same way as in Theorem \ref{DAthrm},  we get
$$
f(x_k) -f(x^*)   =\mathcal{O} \Big(\frac{1}{k}\Big).
$$
\end{proof}

\subsection{Dual Averaging in a distributed setting}
The distributed version of FW-SDA is given in Algorithm 5.  We assume the same setup as in Section 2.3, namely Assumption 1  holds true for the network.
\begin{algorithm}[H]
\textbf{Input :} Objective Function  $f=\frac{1}{m}\sum_{i=1}^m f^i$ ; parameters $\alpha_k$.\\

\textbf{Initialize :} $x^i_0 \in \Omega$, $w^i_0 = \text{argmin}_{x \in \Omega} \langle \nabla f^i(x_0), x \rangle $ and $g^i_0 =0$ for all $i$.\\

\textbf{For} $k =1,2...$ \textbf{do} : At each node $i$,\\ 
\begin{enumerate}

\item  Set $$z^i_{k-1} = (1-\alpha_k)x^i_{k-1} + \alpha_k w^i_{k-1}.$$

\item Compute  the gradient $\nabla f^{i}(z^i_{k-1})$.\\

\item Update the weighted average gradient using the neighbour estimates :
$$g^i_k = \frac{1}{\sum_{p=1}^k \beta_p}\sum_{j=1}^{m} q_{ij } \Big\{ (\sum_{p=1}^{k-1} \beta_p)   g^j_{k-1} + \beta_k \nabla f^{j}(z^j_{k-1})\Big\}. $$

\item Call the LO oracle to compute :
$$w^i_k \in  \text{argmin}_{w\in \Omega} \, \langle g^i_k , w \rangle, $$

\item Set $$ x^i_{k} =(1-\alpha_k) \sum_{j=1}^m q_{ij }x^j_{k-1} + \alpha_k w^i_{k}.$$ 

\end{enumerate}

\textbf{end for }

\caption{Distributed dual averaging for Frank-Wolfe}
\end{algorithm}

We use the stacked vector notation in this section. Specifically for any bold faced variable $\x$, we have $\x := [x^1,....,x^m]$. Also, we let $\f(\x) := [f^1(x^1),....,f^m(x^m)]$.  We let $\varotimes$ denote the
Kronecker product between two matrices . Let $ \bar{\beta}_k = \sum_{i=1}^{k}\beta_i$. . Then, the main steps of Algorithm 5 can be written in a stacked vector notation as follows :
\begin{equation}\label{dist-1}
\z_{k-1} =  (1-\alpha_k) \x_{k-1} + \alpha_k \w_{k-1}.
\end{equation}
\begin{equation}\label{dist-2}
\g_{k+1} = \frac{1}{ \bar{\beta}_k}(Q \varotimes I_d )\{  \bar{\beta}_{k-1} \g_{k}  + \beta_k \nabla \f (\z_{k-1}) \}.
\end{equation}
\begin{equation}\label{dist-3}
\x_{k} =  (Q \varotimes I_d ) (1-\alpha_k )\x_{k-1} + \alpha_k \w_k.
\end{equation}
Let $\LL_{f}  (\z_{k-1}; \w_{k} ) := [L_{f^1}  (z^1_{k-1}; w^1_{k} ) ,...,L_{f^m}  (z^m_{k-1}; w^m_{k} ) ]$, where $L_{f}(\cdot,\cdot)$ is as in (\ref{linear}). We define the distributed average linear model as,
\begin{equation} \label{p1}
\overline{ \PH}_k (\w) := \frac{1}{ \bar{\beta}_k} \sum_{t=1}^{k} Q^{k-t+1} \beta_k  \{\LL_{f}  (\z_{k-1}; \w)  \} 
=\frac{1}{ \bar{\beta}_k} \sum_{t=1}^{k} Q^{k-t+1}\beta_k \{ \f(\z_{t-1}) + \langle  \nabla \f( \z_{t-1}), \w- \z_{t-1}\rangle  \}
\end{equation}
where $\langle  \nabla \f( \z_t), \w- \z_t\rangle := [\langle  \nabla f^1( z_t^1), w^1 -  z^1_t \rangle  ,...,\langle  \nabla f^m( z^m_t), w^m- z^m_t\rangle   ]$. Also, we can write (\ref{dist-2}) as
\begin{align*}
\g_{k+1} &= (Q \varotimes I_d ) \frac{1}{ \bar{\beta}_k} \{ \bar{\beta}_{k-1} \g_{k}  + \beta_k\nabla \f (\z_{k-1})\}, \\
\g_{k+1} &=  \frac{1}{ \bar{\beta}_k} \Big\{  \bar{\beta}_{k-1} (Q^2 \varotimes I_d )\g_{k-1} + \beta_{k-1} (Q^2 \varotimes I_d ) \nabla \f (z_{k-2})   +  \beta_k (Q \varotimes I_d ) \nabla \f (z_{k-1})\Big\}.
\end{align*}
Iterating the above equation we get,
$$ \g_{k+1} = \frac{1}{ \bar{\beta}_k} \Big\{\beta_0(Q^{k+1} \varotimes I_d )\g_{0} + \sum_{t=1}^{k}\beta_t(Q^{k-t+1} \varotimes I_d ) \nabla \f (\z_{t-1}).\Big\} $$
We can assume without loss of generality that $g^i_0 = 0$ for all $i$. Comparing the above equation with (\ref{p1}), we have from the definition of $\w_k$ (Step (d)) \footnote{The equality is interpreted component-wise, i.e. $w_k^i  \in \text{argmin}_{\w}  \overline{\PH}^i_k (w) $ }, 
\begin{equation}\label{min-dist}
\w_k \in \text{argmin}_{\w} \overline{\PH}_k (\w).
\end{equation}

\begin{thm} \label{DAthrmdist}

With $\alpha_k = \frac{2}{k+1} $ and any $\beta_k$ such that $
\alpha_k  = \frac{\beta_k}{\sum_{i=1}^{k}\beta_i},$, we have
$$  f(\overline{x}_k) - f(x^*) = \mathcal{O}(\frac{1}{k}),   $$
where 
$$\bar{x}_k = (Q^* \varotimes  I_d)\, \x_k $$
 and $Q^*= \frac{\boldsymbol{11}^T}{m}$.
\end{thm}

\begin{proof} Before proceeding with the main proof, we prove the following claim :\\

\noindent \textit{Claim :} $\|\x_k - (Q^* \varotimes  I_d) \x_k\| = \mathcal{O} (\alpha_k)$ for all $k \geq 1$.\\
%
%(ii)  $\| \overline{ \PH}_k (\w)  -  \PH_k (\w) \| = \mathcal{O}(\frac{1}{k})  $  for all $k\geq 1$ and any $\w \in \Omega$.\\

\noindent \textit{Proof :} We prove this by induction. Suppose that $$\|\x_k -(Q^* \varotimes  I_d) \x_k\| \leq \alpha_k.$$ We have from (\ref{dist-3}),
\begin{align*}
\x_{k+1} &= (Q \varotimes  I_d)(1-\alpha_k) \x_k+ \alpha_k \w_{k+1}\\
\Longrightarrow    (Q^* \varotimes  I_d) \x_{k+1} &= (Q ^*\varotimes  I_d)(1-\alpha_k) \x_k+ \alpha_k(Q^* \varotimes  I_d) \w_{k+1},
\end{align*}
where we have used the fact that $ (Q^* \varotimes  I_d)\times (Q \varotimes  I_d)  =  (Q^* \varotimes  I_d)$. Then we have,
\begin{multline*}
 \| \x_{k+1} - (Q^* \varotimes  I_d) \x_{k+1} \|  =  \|(Q \varotimes  I_d)(1-\alpha_k) \x_k+ \alpha_k \w_{k+1}\\
  - (Q ^*\varotimes  I_d)(1-\alpha_k) \x_k - \alpha_k(Q^* \varotimes  I_d) \w_{k+1}\|,
\end{multline*}
so,
\begin{multline*}
 \| \x_{k+1} - (Q^* \varotimes  I_d) \x_{k+1} \| \leq (1-\alpha_k) \|\{(Q- Q^*) \varotimes  I_d\} \x_k\| \\
  + \alpha_k  \|\w_{k+1} - (Q^* \varotimes  I_d) \w_{k+1}\|.
\end{multline*}
Noting that $\{(Q-Q^*)\varotimes I_d \} (Q^*\varotimes I_d) x_k =0$, we get 
\begin{multline*}
 \| \x_{k+1} - (Q^* \varotimes  I_d) \x_{k+1} \| \leq (1-\alpha_k) \|(Q \varotimes  I_d) - (Q^* \varotimes  I_d) \| \times \\ \|\x_k - (Q^*\varotimes I_d) \x_k \|
  + \alpha_k  \| \w_{k+1} - (Q^* \varotimes  I_d) \w_{k+1}\|.
\end{multline*}
Since $1-\alpha_k \leq 1 $ and $\|  (Q \varotimes  I_d) - (Q^* \varotimes  I_d) \| \leq \theta \beta $ (from Lemma 1), we have
\begin{multline*} 
 \| \x_{k+1} - (Q^* \varotimes  I_d) \x_{k+1} \| \leq (1-\alpha_k)\theta \beta \|\x_k - (Q^*\varotimes I_d) \x_k \|  + \alpha_k \|\w_{k+1} - (Q^* \varotimes  I_d) \w_{k+1}\|. 
\end{multline*}
The first term is $\mathcal{O}(\alpha_k)$ by the induction step and so is the second one, since $\w_{k+1}$ is bounded. The claim follows.
\\ 
%
%(ii)  Set $R_\Omega := \sup_{\w,\z} \| \f(\z) + \langle  \nabla \f( \z), \w- \z\rangle  \}\|$. Note that $R_\Omega<\infty$. Subtracting (\ref{p2}) from (\ref{p1}) we get,
%\begin{align*}
%\| \overline{ \PH}_k (\w) - \Phi_k (\w)\| &\leq  \frac{1}{k} \sum_{t=1}^{k} \| \big( Q^{k-t+1}-Q^*  \big) \varotimes I_d    \|\Big( \|\f(\z_{t-1})+  \nabla \f( z_{t-1}) \langle \w-\z_t\rangle \|   \Big),  \\
%& \leq  \frac{R_\Omega}{k} \sum_{t=1}^{k}  \theta \beta^{k-t+1}    \qquad (\text{Lemma 1})\\
%&\leq \frac{R_\Omega \theta \beta}{(1-\beta)k} =\mathcal{O}(\frac{1}{k}).
%\end{align*}

We continue with the proof of the theorem. From (\ref{smooth}), we have for all $i$ :
\begin{align*}
f^i(x^i_k) &\leq L_{f^i} (z^i_{k-1}; x^i_k ) + \frac{L}{2} \|x^i_k -z^i_{k-1} \| ^2,\\ 
&= (1 - \alpha_k) L_{f^i} (z^i_{k-1}; x^i_{k-1} ) + \alpha_k L_{f^i}  (z^i_{k-1}; w^i_{k} ) + \frac{L}{2} \|x^i_k -z^i_{k-1} \| ^2.
\end{align*} 
The convexity of $f^i(\cdot)$ gives
\begin{equation*}
f^i(x^i_k)  \leq (1 - \alpha_k) f(x^i_{k-1})  + \alpha_k L_{f^i}  (z^i_{k-1}; w^i_{k} ) + \frac{L}{2} \|x^i_k -z^i_{k-1} \| ^2.
\end{equation*}
Summing over $j$, we get
\begin{equation*}
\sum_{j=1}^{m} q_{ij} f^j(x^j_k)  \leq (1 - \alpha_k) \sum_{j=1}^{m}q_{ij} \big\{ f^j(x^j_{k-1})   + \alpha_k L_{f^j}  (z^j_{k-1}; w^j_{k} ) +  \frac{L}{2} \|x^i_k -z^i_{k-1} \| ^2\big\},
\end{equation*}
so that in vector notation,
\begin{equation}\label{dist-4}
Q\f(\x_k)  \leq (1 - \alpha_k) Q \f( \x_{k-1})  + \alpha_k Q \LL_{f}  (\z_{k-1}; \w_{k} )  + \frac{L}{2}   \| \x_k -\z_{k-1} \| ^2 \textbf{1},
\end{equation}
where the inequality is interpreted component-wise. Subtracting (\ref{dist-1}) from (\ref{dist-3}), we get
$$
\| \x_k-\z_{k-1} \| =  (1-\alpha_k )\|\x_{k-1} - (Q\varotimes I_d) \x_{k-1}\|+ \alpha_k \|\w_k-\w_{k-1} \|.
$$
Using Claim (i), we have $\| \x_k- \z_{k-1} \| = \mathcal{O}(\alpha_k)$, so that $\| \x_k-\z_{k-1} \| \leq \sqrt{C_{\Omega}'} \alpha_k$ for some constant $C_{\Omega}'$. Using this in (\ref{dist-4}),
\begin{equation}\label{dist-5}
Q \f(\x_k)  \leq (1 - \alpha_k) Q \f( \x_{k-1})  + \alpha_k Q \LL_{f}  (\z_{k-1}; \w_{k} )  + \frac{L C_{\Omega}' \alpha^2_k }{2} \textbf{1}.
\end{equation} 
We also have,
 $$  \bar{\beta}_k \overline{\PH}_k (\w_k)   =  \beta_kQ \LL_{f} (\z_{k-1} ; \w_k ) + \bar{\beta}_{k-1} \overline{\PH}_{k-1} (\w_k) . $$
From (\ref{min-dist}), $ \w_{k-1}  \in \text{argmin}_w \overline{\PH}_{k-1} (w) $, so that $ \overline{\PH}_{k-1}( \w_{k-1}) \leq   \overline{\PH}_{k-1} (\w_k)$. This gives 
\begin{align*}
\bar{\beta}_k  \overline{\PH}_k (\w_k)  &\geq  \beta_k Q \LL_{f} (\z_{k-1} ; \w_k ) + \bar{\beta}_{k-1} \overline{ \PH}_{k-1} (\w_{k-1}) ,\\  
 \Longrightarrow \beta_k Q  \LL_{f} (\z_{k-1} ; \w_k )&\leq   \bar{\beta}_{k} \overline{\PH}_k (\w_k) -  \bar{\beta}_{k-1}\overline{ \PH}_{k-1} ( \w_{k-1}) .
\end{align*} 
Proceeding in the same way as Theorem 1, 
\begin{align*}
\alpha_k \beta^{-1}_k  [ \bar{\beta}_{k} \overline{ \PH}_k (\w_k) -  \bar{\beta}_{k-1}\overline{ \PH}_{k-1} (\w_k)] &=  \overline{ \PH}_k (\w_k)- \alpha_k \beta^{-1}_k \bar{\beta}_{k-1}\   \overline{ \PH}_{k-1} (\w_{k-1})],\\
&= \overline{ \PH}_k (\w_k) - \alpha_k \beta_k^{-1}(\bar{\beta}_k -\bar{\beta}_{k-1}) \overline{ \PH}_{k-1} (\w_{k-1}),\\
&=  \overline{ \PH}_k (\w_k) -   ( 1-\alpha_k ) \overline{ \PH}_{k-1} (\w_{k-1}).
\end{align*} 
Combing the above facts we have,
\begin{equation*}
\alpha_k Q \LL_{f} (\z_{k-1} ; \w_k ) \leq  \overline{ \PH}_k (\w_k) -  (1-\alpha_k ) \overline{ \PH}_{k-1} (\w_{k-1}).
\end{equation*}
Using this in (\ref{dist-5}) we get,
\begin{equation*}
Q\f(\x_k)  \leq (1 - \alpha_k) Q \f(\x_{k-1})  + \big(  \overline{ \PH}_k (\w_k) -  (1-\alpha_k ) \overline{ \PH}_{k-1} (\w_{k-1}) \big)  + \frac{L C_{\Omega}' \alpha^2_k }{2} \textbf{1} .
\end{equation*}
Rearranging, we get
\begin{equation}\label{dist-6}
Q \f(\x_k) -\overline{ \PH}_k (\w_k)   \leq (1 - \alpha_k) \{ Q \f(\x_{k-1})  - \overline{ \PH}_{k-1} (\w_{k-1}) \} + \frac{L}{2}\alpha^2_k  C_{\Omega}^2 \cdot \textbf{1} .
\end{equation}
Set $ \nu_k =  Q\f(\x_k) -\overline{ \PH}_k (\w_k)   $ and  define 
 \begin{equation*}
   \Gamma_k :=
    \begin{cases}
      (1-\alpha_k)  \Gamma_{k-1}, & \text{if}\ k\geq 2 \\
      1, &  \text{if}\ k=1.
    \end{cases}
  \end{equation*}
We can write (\ref{dist-6}) with the above notation as\footnote{ interpreted component-wise} :
\begin{align*}
\nu_k &\leq (1 - \alpha_k) \nu_{k-1} + + \frac{L}{2}\alpha^2_k  C_{\Omega}^2, \\
\Longrightarrow \frac{\nu_k }{\Gamma_k} &\leq  \frac{ \nu_{k-1}}{\Gamma_{k-1}} + \frac{L}{2}\frac{\alpha^2_k}{\Gamma_k}  C_{\Omega}^2  \,\,\,\,\, (\text{since } \Gamma_k ).
\end{align*}
A simple calculation gives :
$$
\Gamma_k = \frac{2}{k(k+1)} \,\,\,\text{       and        }\,\,\,\,  \frac{\alpha_k^2}{\Gamma_k}\leq 2.
$$
To conclude the proof we note that,
\begin{align*}
\nu_k &\leq \Gamma_k \Big\{ \frac{\nu_1}{\Gamma_1} + \sum_{i=1}^{k} \big( \frac{L}{2}  C_{\Omega}^2 + C_{\Omega} \big)  \Big( \frac{\alpha^2_i}{\Gamma_i} \Big) \Big\},\\
\nu_k &\leq \frac{2}{k(k+1)}\Big\{ \frac{\nu_1}{\Gamma_1} + k(L C_{\Omega}^2 + 2C_{\Omega})\Big\},\\
& = \mathcal{O}(\frac{1}{k}).
\end{align*}
Thus we have,
$$
Q \f(\x_k) - \overline{ \PH}_k (\w_k) \leq \mathcal{O}(\frac{1}{k}).
$$
Since $Q^*\times Q =Q^*$, we have
\begin{equation}\label{aux-10}
Q^* \f(\x_k) - Q^*  \overline{ \PH}_k (\w_k) \leq \mathcal{O}(\frac{1}{k}).
\end{equation}
Let $\x^*:= [x^*,...,x^*]$. We have, $$   \overline{ \PH}_k (\w_k)  \leq \overline{ \PH}_k (\x^*) \leq \frac{1}{\bar{\beta_k}} \sum_{t=1}^{k} \beta_t Q^{k-t+1}  \f(\x^*)$$
so that,
 $$  Q^* \overline{ \PH}_k (\w_k) \leq Q^*  \f(\x^*).$$
Using the above inequality in (\ref{aux-10}), we get,
$$
Q^* \f(\x_k) - Q^*\f(\x^*) \leq  Q^* \f(\x_k)  -  Q^* \overline{ \PH}_k (\w_k) \leq \mathcal{O}(\frac{1}{k}).
$$
Using any row of the above inequality we get 
$$
\frac{1}{m}\sum_{i=1}^mf^i( x^i_k) -\frac{1}{m}\sum_{i=1}^mf^i( x^*)\leq \mathcal{O}(\frac{1}{k}).
$$
Then using claim and the $L$-smoothness of $f(\cdot)$,
\begin{multline*}
|f(\overline{\x}_k)-f({x}^*) |= | f(\overline{\x}_k)- \frac{1}{m}\sum_{i=1}^m f^i(x^i_k) |+|\frac{1}{m} \sum_{i=1}^m f^i(x^i_k)  -f({x}^*)|  \leq  \mathcal{O}(\alpha_k) +\mathcal{O}(\frac{1}{k}) 
=  \mathcal{O}(\frac{1}{k}),
\end{multline*}
which concludes the proof.
\end{proof}

\bibliographystyle{IEEEtran}
\bibliography{biblo}

\end{document}